\documentclass[10pt]{article}
\usepackage{amsmath}
\usepackage{english}
\usepackage{amsfonts}
\usepackage{amssymb}
\usepackage{a4wide}

\newtheorem{theoreme}{Th\'eor\`eme}[section]

\newtheorem{proposition}[theoreme]{Proposition}
\newtheorem{corollaire}[theoreme]{Corollaire}
\newtheorem{definition}[theoreme]{D\'efinition\rm}
\newtheorem{remarque}{\it Remarque}

 \def\R{\mathbb{R}}
\def\QQQ{\mathbb{Q}}

\def\Z{\mathbb{Z}}

\def\Rplusetoile{\R _+ ^\ast}

\def\ssi{si, et seulement si, }

\newcommand\al{\alpha}
\newcommand\strois{{\mathfrak S}_3}
\newcommand\scinq{{\mathfrak S}_5}

\newcommand{\dcinq}{{\cal D}_5}

\newcommand{\croix}{\times}

\newcommand{\dd}{{\rm d} }    

\newcommand{\psd}{\rtimes}

\newcommand{\zdeuxz}{\Z / 2 \Z}

\newcommand{\xsoul}{\underline{x}}

\newcommand{\AAA}{A}
\newcommand{\BBB}{B}
\newcommand{\CCC}{C}

\newcommand{\zeron}{\{0,\ldots,n\}}

\newcommand{\unn}{\{1,\ldots,n\}}

\newcommand{\unnmd}{\{1,\ldots,n-2\}}
\newcommand{\unnmdeux}{\{1,\ldots,n-2\}}

\newcommand{\deuxn}{\{2,\ldots,n\}}

\newcommand{\deuxnmd}{\{2,\ldots,n-2\}}

\newcommand{\troisn}{\{3,\ldots,n\}}

\newcommand{\troisnmu}{\{3,\ldots,n-1\}}
\newcommand{\zupn}{{[0,1]^n}}

\newcommand{\de}{\delta}

\newcommand{\psii}{\psi}
\newcommand{\siigma}{\sigma}
\newcommand{\psiexp}{{{\psi}}}
\newcommand{\phiexp}{{{\phi}}}

\newcommand{\sigmaexp}{{{\siigma}}}

\newcommand{\congru}{\equiv}
\newcommand{\noncongru}{\not\equiv}

\newcommand{\df}{D}
\newcommand{\JJJ}{L}
\newcommand{\KKK}{J}
\newcommand{\LLL}{K}
\newcommand{\ro}{\varrho}

\newcommand{\XXX}{y}
\newcommand{\Li}{{\rm Li}}

\newcommand{\zuptrois}{{[0,1]^3}}

\newcommand{\XXXpr}{y'_}

\newcommand{\ddeux}{V}
\newcommand{\petip}{\underline{p}}
\newcommand{\gdp}{\underline{P}}

  \newcommand{\GL}{{\rm GL}}
 \newcommand{\eee}{{\cal E}}
\newcommand{\phii}{\phi}
   \newcommand{\eps}{\varepsilon}

\title{Formes lin\'eaires en polyz\^etas et int\'egrales multiples}
\author{St\'ephane Fischler} 
\date{7 F\'evrier 2002}

\begin{document}
\maketitle

 \section{Introduction} \label{secintro}

Apr\`es la d\'emonstration de l'irrationalit\'e de $\zeta(3)$ par Ap\'ery 
\cite{Apery}, plusieurs variantes ont
\'et\'e propos\'ees, parmi lesquelles celles de Beukers \cite{Beukers} et
 Sorokin \cite{Sorokin}. Dans ces deux preuves, les formes lin\'eaires
  en 1 et $\zeta(3)$ sont \'ecrites comme des int\'egrales
triples : 
$$B(N) = \int_{\zuptrois} \frac{x^N (1-x)^N y^N (1-y)^N 
z^N (1-z)^N}{(1-z(1-y(1-x)))^{N+1}}
\dd x \dd y \dd z$$
 pour Beukers, et
$$S(N)=\int_{\zuptrois} \frac{x^N (1-x)^N y^N (1-y)^N z^N (1-z)^N}{(1-xy)^{N+1}
 (1-xyz)^{N+1}} \dd x \dd y \dd z$$
  pour Sorokin. Or ces formes lin\'eaires
 co\"{\i}ncident. Si on croit \`a la philosophie des p\'eriodes 
 \cite{KontsevichZagier}, l'\'egalit\'e de ces int\'egrales doit pouvoir se 
 d\'emontrer par une suite de changements de variables et d'applications
 des r\`egles d'additivit\'e (par rapport \`a l'int\'egrande ou au domaine)
 et du th\'eor\`eme de Stokes. En l'occurrence, un seul changement de variables
 suffit (voir le th\'eor\`eme \ref{theocdv}). 
 
Les int\'egrales $B(N)$ ont \'et\'e g\'en\'eralis\'ees par 
Vasilyev \cite{Vasilyevancien}, qui   pose  $\de_k(x_1,\ldots,
x_n) = 1-x_k \de_{k-1} (x_1,\ldots,
x_n)$ pour $n \geq 2$ et $k \in \unn$ avec $\de_0 = 1$, et   consid\`ere 
\begin{equation} \label{intVasilyev}
\int_{\zupn} 
 \frac{\prod_{k=1} ^n x_k ^N (1-x_k)^N}{\de_n(x_1,\ldots,x_n) ^{N+1}}
\dd x_1 \ldots \dd x_n.
\end{equation}
Il d\'emontre que pour $N=0$ cette int\'egrale est un multiple rationnel de
$\zeta(n)$, et \cite{Vasilyev}
que pour $n=5$ (respectivement $n=4$)
et $N$ quelconque c'est une forme lin\'eaire 
en 1, $\zeta(3)$ et $\zeta(5)$ (respectivement 
1, $\zeta(2)$ et $\zeta(4)$). Il conjecture que pour tous $n$ et $N$ on 
obtient une forme 
lin\'eaire en 1 et les valeurs de $\zeta$ aux entiers compris entre 2 et $n$
ayant la m\^eme parit\'e que $n$.

D'autre part, les int\'egrales $S(N)$ sont \`a rapprocher
de celles que Sorokin introduit  \cite{Sorokinpi}  
 en relation avec $\zeta(2,2,\ldots,2)$,
 quand $n=2r$ est pair :
\begin{equation} \label{intSorokinpi}
\int_{\zupn} 
 \frac{(\XXX_1\XXX_2)^{rN+(r-1)} (\XXX_3\XXX_4)^{(r-1)N+(r-2)}\ldots
 (\XXX_{n-1}\XXX_n)^N
  \prod_{k=1} ^n  (1-\XXX_k)^N}{\prod_{{\tiny k \in \deuxn {\rm pair}}}
  (1-\XXX_1\XXX_2\ldots \XXX_k)^{N+1}}
\dd \XXX_1 \ldots \dd \XXX_n.
\end{equation}

Dans la section \ref{secun} ci-dessous, on d\'efinit deux familles d'int\'egrales, qui
g\'en\'eralisent $B(N)$ et (\ref{intVasilyev}) d'une part,
$S(N)$ et  (\ref{intSorokinpi}) d'autre part, et on montre 
que ces deux familles se correspondent par un changement de variables.
On montre en outre qu'un groupe agit sur ces int\'egrales, de mani\`ere
 analogue
\`a ce que Rhin et Viola consid\`erent (\cite{RV2}, \cite{RV3}) 
pour obtenir les meilleures mesures d'irrationalit\'e  connues pour $\zeta(2)$
et $\zeta(3)$. \`A la section \ref{secdeux}, on d\'efinit une autre famille
d'int\'egrales $n$-uples, qui g\'en\'eralise (\ref{intVasilyev}) et
sur laquelle agit aussi un groupe ; on retrouve
alors dans les cas particuliers   $n=2$
et $n=3$ les groupes obtenus par  Rhin et Viola.

On adopte la d\'efinition suivante :
\begin{definition} On dit qu'une famille de nombres $I(\petip) \in \Rplusetoile
\cup \{ \infty\}$, param\'etr\'es par $\petip \in \Z^s$, admet pour groupe de 
Rhin-Viola un sous-groupe $G$ de $\GL_s(\Z)$ si $I(g \petip) / I(\petip)$ est 
un rationnel non nul pour tous $g \in G $ et $\petip \in \Z^s$ (avec $\infty
 / \infty = 1$).
\end{definition}

Dans ce texte, on consid\`ere des familles d'int\'egrales 
$n$-uples dont les valeurs
(lorsqu'elles sont finies) sont conjecturalement des formes lin\'eaires sur 
$\QQQ$ en les polyz\^etas de poids au plus $n$. On peut esp\'erer que, pour
de bons choix des exposants $\petip$, ces int\'egrales soient assez
petites et qu'on ait un  
contr\^ole sur le d\'enominateur des coefficients de la forme lin\'eaire. Alors
l'\'etude de la valuation $p$-adique des nombres 
$I(g \petip) / I(\petip)$, quand $g$ parcourt $G$ et $p$ un certain ensemble de
nombres premiers, peut permettre d'am\'eliorer
ce contr\^ole du d\'enominateur, donc de raffiner des mesures
d'irrationalit\'e ou d'ind\'ependance lin\'eaire de certains polyz\^etas.

Davantage de d\'etails, en particulier sur l'action des groupes
de Rhin-Viola, seront donn\'es dans \cite{SFCaen}.
 
\section{Une g\'en\'eralisation commune des int\'egrales de Vasilyev et de Sorokin} \label{secun}


Dans tout ce texte, $n$ d\'esigne  un entier sup\'erieur ou \'egal \`a 2. 

Pour $k \in \zeron$  et $\xsoul=(x_1,\ldots,x_n) \in \zupn$ on pose
$\df_k(\xsoul) = \sum_{j=0} ^k (-1)^j x_n x_{n-1}   \ldots   x_{n-j+1}$.
On a alors $\df_0(\xsoul)=1$, $\df_1(\xsoul)=1-x_n$,
 $\df_2(\xsoul)=1-x_n(1-x_{n-1})$ et $\df_n(\xsoul)=\de_n(\xsoul)$.
 \`A tout $\petip = 
 (a_1,\ldots,a_n,b_1,\ldots,b_n,c_2,\ldots,c_n) \in \Z^{3n-1}$  
   on associe l'int\'egrale (\'eventuellement infinie)
$$\KKK(\petip) = \int_{[0,1]^n}
 \frac{\prod_{k=1} ^n x_k ^{a_k} (1-x_k)^{b_k}}{\prod_{k=2} ^n 
 \df_k(\xsoul) ^{c_k}}
 \frac{\prod_{k \in \deuxnmd \mbox{{\small pair }}} \df_k(\xsoul)}{\prod_{k 
\in \troisnmu \mbox{{\small impair }}} \df_k(\xsoul)}  
\frac{\dd x_1  \ldots  \dd x_n}{\df_n(\xsoul)}. 
$$

Par ailleurs, \`a tout $\gdp =   (\AAA_1,\ldots,\AAA_n,\BBB_1,\ldots,\BBB_n,
\CCC_2,\ldots,\CCC_n) \in \Z^{3n-1}$  
     on associe 
$$\LLL(\gdp) = \int_{[0,1]^n}
 \frac{\prod_{k=1} ^n \XXX_k ^{\AAA_k} (1-\XXX_k)^{\BBB_k}}{\prod_{k=2} ^n 
 (1-\XXX_1\XXX_2\ldots \XXX_k) ^{\CCC_k+1}} 
    \dd \XXX_1   \ldots   \dd \XXX_n. 
$$
Cette int\'egrale est finie 
\ssi on a  
$ \sum_{k=2} ^n (\CCC_k - \BBB_k)^+  \leq \BBB_1$ et $\AAA_k \geq 0$, 
$\BBB_k \geq 0$ pour tout $k$.

\bigskip

\begin{theoreme} \label{theocdv}
Pour tout $\petip \in \Z^{3n-1}$ on a $\KKK(\petip)=
\LLL(\gdp)$, o\`u $\gdp$ est donn\'e en fonction de $\petip$ par :
$$\begin{array}{l}
\AAA_k = a_{n+1-k}\mbox{ pour } 1 \leq k \leq n \mbox{ ,}\\
\BBB_k =b_{n+1-k} \mbox{ pour } 2 \leq k \leq n \mbox{ et }
\BBB_1 = a_{n-1} + b_n - c_2 - c_3 - \ldots -c_n, \\
\CCC_k =a_{n+1-k} + b_{n+1-k} - c_k - c_{k+1} - \ldots -c_n
        \mbox{ pour tout } k \in \deuxn \mbox{ pair,}\\
\CCC_k =   c_k + c_{k+1} + \ldots +c_n - a_{n-k}
        \mbox{ pour tout } k \in \troisn \mbox{ impair, avec la convention
	$a_0 = 0$.}
\end{array}$$
\end{theoreme}

Ce r\'esultat provient du changement de variables d\'efini par
$x_k  =  \XXX_{n+1-k}$ pour $k \congru n \mod 2$ et $x_k  = 
 \frac{(1-\XXX_1 \ldots \XXX_{n-k})\XXX_{n+1-k}}{1-\XXX_1 
 \ldots \XXX_{n+1-k}}$ pour
 $k \noncongru n \mod 2$. On peut bien s\^ur l'inverser, ce qui correspond 
 \`a exprimer $\petip$ en fonction de $\gdp$ dans le th\'eor\`eme
 \ref{theocdv}. 

\begin{corollaire} Avec les notations de l'introduction, on a $B(N)=S(N)$
pour tout $N \geq 0$, et l'int\'egrale (\ref{intVasilyev}) est \'egale \`a
$\left\{ \begin{array}{l}
\int_{\zupn} 
 \frac{\prod_{k=1} ^n \XXX_k ^N (1-\XXX_k)^N}{\prod_{k \in \deuxn \mbox{{\small
 pair}}} (1-\XXX_1 \ldots \XXX_k) ^{N+1}}
\dd \XXX_1 \ldots \dd \XXX_n \mbox{ si } n  \mbox{ est pair,} \\
\int_{\zupn} 
 \frac{\prod_{k=1} ^n \XXX_k ^N (1-\XXX_k)^N}{(1-\XXX_1 \ldots \XXX_n) ^{N+1}
 \prod_{k \in \deuxn \mbox{{\small
 pair}}} (1-\XXX_1 \ldots \XXX_k) ^{N+1}}
\dd \XXX_1 \ldots \dd \XXX_n  \mbox{ si } n  \mbox{ est impair.}
\end{array} \right. $
\end{corollaire}

\medskip

L'int\'er\^et des  int\'egrales  $\LLL(\gdp)$ (donc du th\'eor\`eme
 \ref{theocdv}) est qu'elles se d\'eveloppent ``naturellement'' en 
 s\'eries multiples, dont on peut esp\'erer d\'emontrer que ce sont des
   formes lin\'eaires sur $\QQQ$
en les polyz\^etas de poids au plus $n$. Par exemple, pour $N=0$,
(\ref{intVasilyev}) vaut 
$$ \sum_{l_1\geq l_2 \geq \ldots
\geq l_{(n-1)/2} \geq l_{(n+1)/2}  \geq 1} \frac{1}{l_1 ^2 l_2 ^2  \ldots 
  l_{(n-1)/2}^2  l_{(n+1)/2} }$$  si $n$ est impair, et
  $$ \sum_{l_1\geq l_2 \geq \ldots
\geq l_{n/2} \geq 1} \frac{1}{l_1 ^2 l_2 ^2    \ldots l_{n/2} ^2}$$
 si $n$ est pair. 
Le r\'esultat de Vasilyev selon lequel cette int\'egrale \'egale  
$(-1)^{n-1} \Li_n((-1)^{n-1})$ se 
 ram\`ene ainsi \`a  une identit\'e lin\'eaire
entre
polyz\^etas. Or il existe des outils combinatoires pour d\'emontrer de telles
identit\'es (voir \cite{MiW}).

\begin{proposition} \label{propgrpL}
Si $n \geq 3$, la famille $(\LLL(\gdp))$ admet un groupe de Rhin-Viola d'ordre
32, isomorphe \`a $( \ddeux \croix \ddeux) \psd \zdeuxz$, o\`u 
$\ddeux = \zdeuxz \croix \zdeuxz $ est 
le groupe de Klein et $\zdeuxz$ agit en permutant les facteurs de
$\ddeux \croix \ddeux$.
 Cela provient d'analogues des transformations $\sigma$
et $\chi$ de \cite{RV3}, et du changement de variables d\'efini par
$ \XXXpr{k} =\frac{1-\XXX_1\ldots \XXX_{n-k+1}}{1-\XXX_1\ldots \XXX_{n-k+2}}$
pour tout $k \in \unn$ (avec la convention $\XXX_{n+1}=0$).
\end{proposition}

\section{Une g\'en\'eralisation du groupe de Rhin-Viola} \label{secdeux}

On reprend la notation $\de_k$ utilis\'ee dans l'introduction.
\`A tout $\petip = (a_1,\ldots,a_n,b_1,\ldots,b_n,c_2,\ldots,c_n) \in
 \Z^{3n-1}$  on associe
$$\JJJ(\petip) = \int_{[0,1]^n}
 \frac{\prod_{k=1} ^n x_k ^{a_k} (1-x_k)^{b_k}}{\prod_{k=2} ^n
  \de_k(\xsoul) ^{c_k}}   \frac{\dd x_1  \ldots   \dd x_n}{\de_n(\xsoul)}.
$$
On pose $\ro_n = c_n-b_n$ et $\ro_{n-1} =  c_{n-1} -1 -b_{n-1}$,
puis $\ro_k = \ro ^+ _{k+2} + c_k -1 -b_k$ pour tout $k \in \unnmdeux$ 
en notant $\al ^+ = \max(\al,0)$ et avec la convention $c_1 = 1$. 
Alors l'int\'egrale $\JJJ(\petip)$   est finie \ssi on a $a_k \geq 0$,  
$b_k \geq 0$ et  $\ro_k \leq a_{k-1}$
pour tout $k \in \unn$, avec la convention
$a_0 = 0$.

\bigskip

Notons $\siigma$ l'automorphisme de $\Z^{3n-1}$ qui \'echange $a_1$ et $b_2$,
ainsi que $a_2$ et $b_1$, en fixant les autres coordonn\'ees. Notons 
$\psii$ celui qui \`a $\petip$ associe $\petip ' =(a'_1,\ldots,a'_n,
b'_1,\ldots,b'_n,c'_2,\ldots,c'_n)$ d\'efini par :
$$\begin{array}{l}
a'_k = a_{n+1-k} \mbox{ pour } 1 \leq k \leq n \mbox{ , }
b'_k = b_{n+2-k} \mbox{ pour } 2 \leq k \leq n \mbox{ et }
b'_1 = a_{n-1} + b_n - c_n,\\
c'_k = a_{n+2-k}+b_{n+2-k}+c_{n+1-k}-b_{n+1-k}-a_{n-k} 
		\mbox{ pour } 2 \leq k \leq n-1,\\
c'_n = a_2 + b_2 - b_1.
\end{array}$$
Alors des changements de variables montrent qu'on a $\JJJ(\petip)=
\JJJ(\siigma(\petip))=\JJJ((\psii(\petip))$ pour tout $\petip$.
En outre, notons $\chi$ l'automorphisme de $\Z^{3n-1}$ qui fixe toutes 
les composantes, sauf $a_n$ et $c_n$ qu'il \'echange et $b_n$
qu'il remplace par $a_n + b_n - c_n$ ; il v\'erifie
$ \JJJ(\petip) = \frac{a_{n}! b_{n}!}{ c_n! (a_{n}+b_{n}-c_n) !}
  \JJJ(\chi(\petip))$ pour tout $\petip$.

\begin{proposition} \label{propRV32}
Si $n \geq 3$, la famille $(\JJJ(\petip))$ admet un groupe 
de Rhin-Viola d'ordre
32, isomorphe \`a $( \ddeux \croix \ddeux) \psd \zdeuxz$, qui
est engendr\'e par $\siigma$, $\psii$ et $\chi$.
\end{proposition}
 
\begin{remarque}  Peut-\^etre y a-t-il un lien entre les
int\'egrales  $\JJJ(\petip)$ et $\LLL(\gdp)$ qui permette d'expliquer
le parall\'elisme entre les propositions \ref{propgrpL} et \ref{propRV32} ?
\end{remarque}
 
\bigskip
  
Pour obtenir une structure de groupe plus riche que celle de la proposition
\ref{propRV32}, 
on peut se restreindre aux int\'egrales $\JJJ(\petip)$ telles
 que $c_2 = \ldots = c_{n-1}=0$. Pour conserver l'action de 
 $\siigma$ et $\psii$, on doit imposer en outre $a_1+b_2=a_3+b_3$ 
 si $n=3$, et les 
 relations suivantes si $n \geq 4$ :
$$a_2 = b_1 \mbox{ et } b_n = c_n \mbox{ et }
a_k + b_{k+1} = a_{k+2} + b_{k+2} \mbox{ pour tout }
 k \in \unnmd.$$
On note $\eee$ l'ensemble des $\petip$ v\'erifiant ces relations ; il est stable
par $\siigma$ et $\psii$.
Notons $\phii$
l'automorphisme  de $\eee$ qui stabilise toutes les coordonn\'ees, sauf 
$a_{n-1}$ et $ c_n$ qu'il \'echange, $b_{n-1}$ qu'il remplace par 
$ a_{n-1} + b_{n-1} -c_n$ et $b_n$ qu'il remplace par 
$a_{n-1} + b_n -c_n$. On a alors $\JJJ(\petip) = 
 \frac{a_{n-1}! b_{n-1}!}{ c_n! (a_{n-1}+b_{n-1}-c_n) !}
\JJJ(\phii(\petip))$.

\bigskip

\begin{theoreme} La famille des $\JJJ(\petip)$, pour $\petip \in \eee$, admet
un groupe de Rhin-Viola $G$ engendr\'e par 
$\siigma$, $\psii$ et $\phii$. Plus pr\'ecis\'ement :
\begin{itemize}
\item Pour $n \geq 4$, ce groupe est isomorphe \`a 
$(\strois \croix \strois) \psd \zdeuxz$, donc d'ordre 72 ; il laisse
stable $\frac{\JJJ(\petip)}{ a_{n-1} !   b_{n-1} !   a_2 !  b_3 !}$
si $n \geq 5$, et $\frac{\JJJ(\petip)}{ a_3 !   b_3 !  a_2  ! }$
si $n =4$.
\item Pour $n =3$, le groupe $G$ est isomorphe \`a $H \psd \scinq $, o\`u
$H$ est l'hyperplan $\eps_1 + \ldots + \eps_5 = 0$ de $(\zdeuxz)^5$ ; il laisse
stable $\frac{\JJJ(\petip)}{a_1! a_2! a_3!  b_1! b_2! b_3! (a_2+b_3-c_3)!
(b_1+b_3-c_3)!}$. 
\item Pour $n =2$, le groupe $G$ est isomorphe \`a $\scinq$, et  
laisse
stable $\frac{\JJJ(\petip)}{a_1! a_2!  b_1! b_2! (a_1+b_2-c_2)!}$. 
\end{itemize}
\end{theoreme}

\medskip

Pour $n \in \{2,3\}$, on retrouve exactement les situations consid\'er\'ees
 par Rhin et Viola. Pour $n=2$,
le groupe di\'edral $\dcinq$ de \cite{RV2} est exactement celui engendr\'e
par $\psii$ et $\siigma$. Pour $n=3$, la transformation not\'ee $\vartheta^2$
dans \cite{RV3} est ici $\psii \circ \siigma$. On obtient $\vartheta$
en observant qu'il agit sur $\eee$ comme
$\phiexp (\sigmaexp \phiexp \psiexp \phiexp )^2$. 
 
Pour $n=3$, un
  ph\'enom\`ene
   myst\'erieux se produit dans \cite{RV3} : lorsqu'on impose la relation
    $a_1+b_2=a_3+b_3$, ce qui permet d'avoir une action de groupe, 
 les int\'egrales obtenues sont des formes lin\'eaires en 1 et 
 $\zeta(3)$ seulement : $\zeta(2)$ n'appara\^{\i}t plus.
 On peut se demander si un ph\'enom\`ene analogue
  survient pour $n \geq 4$.

\smallskip

\begin{remarque} Dans tout ce texte, les propri\'et\'es de $\chi$ et $\phii$
proviennent, comme dans \cite{RV2} et \cite{RV3}, de la formule 
$\int_0 ^1 \frac{x^a (1-x)^b}{(1+\beta x)^{c+1}} \dd x
= \frac{a! b! }{c! (a+b-c)!}
\int_0 ^1 \frac{x^c (1-x)^{a+b-c}}{(1+\beta x)^{a+1}} \dd x$. Un analogue de 
cette formule, dans lequel  le d\'enominateur serait de la forme
$(1+\beta x)^{c+1}(1+\beta' x)^{c'+1}$, permettrait d'enrichir les structures
de groupe obtenues ici.
\end{remarque}
 
\hspace{-\parindent}\textsc{Remerciements : }Je remercie 
Pierre Cartier et Tanguy Rivoal, 
dont les questions sont \`a l'origine 
de ce travail, ainsi que Jacky Cresson et Michel Waldschmidt.


\begin{thebibliography}{99}
\selectlanguage{english}
\bibitem{Apery}  Ap\'ery R.,  Irrationalit\'e de $\zeta(2)$ et $\zeta(3)$,
 Ast\'erisque 61 (1979), 11-13. 
\bibitem{Beukers}  Beukers F.,  A note on the irrationality of $\zeta(2)$ 
and $\zeta(3)$,
Bull. London Math. Soc. 11.3  (1979),   268-272. 
\bibitem{SFCaen}   Fischler S.,  Groupes de Rhin-Viola et int\'egrales
multiples, soumis aux Actes des
Rencontres Arithm\'etiques de Caen (juin 2001).
\bibitem{KontsevichZagier} Kontsevich M., Zagier D., Periods, in: Mathematics
Unlimited - 2001 and beyond,   771-808, Springer, 2001. 
\bibitem{RV2}  Rhin G.,  Viola C.,   On a permutation group 
related to $\zeta(2)$, Acta Arith. 77.1 (1996), 23-56.
\bibitem{RV3}  Rhin G.,  Viola C., The group structure for $\zeta(3)$,
 Acta Arith. 97.3 (2001), 269-293.
\bibitem{Sorokinpi}  Sorokin V.N.,   A transcendence measure for $\pi^2$,
Sb. Math. 187:12 (1996), 1819-1852.  
\bibitem{Sorokin}   Sorokin V.N.,  Ap\'ery's theorem, Moscow Univ.
Math. Bull. 53.3 (1998), 48-52.
\bibitem{Vasilyevancien} Vasilyev D.V.,   Some formulas for Riemann 
  zeta-function at  integer points, Moscow Univ. Math. Bull. 51.1 (1996),
  41-43.
\bibitem{Vasilyev}  Vasilyev D.V.,   On small linear forms for the values
of the Riemann zeta-function at odd integers, preprint.
\bibitem{MiW}  Waldschmidt M., Valeurs z\^eta multiples : une introduction,
 J. Th\'eor. Nombres Bordeaux 12.2 (2000), 581-595. 
\end{thebibliography}
\end{document}